\documentclass[11pt]{article}
\usepackage{a4}
\usepackage[T1]{fontenc}
\usepackage[english]{babel}
\usepackage{amsmath}
\usepackage{amssymb}
\usepackage{amsfonts}
\usepackage{mathrsfs}
\usepackage{fig4tex184}
\makeatletter
\@addtoreset{equation}{section} 
\makeatother

\newcommand{\ud}{\mathrm{d}}
\newcommand{\p}{\partial}
\makeatletter
\renewcommand*{\@fnsymbol}[1]{\ensuremath{\ifcase#1\or*\or\ddagger\or\mathsection\or\mathparagraph\or \|\or **\or \dagger\dagger\or \ddagger\ddagger \else\@ctrerr\fi}}
\makeatother
\usepackage{fancyhdr}
\fancypagestyle{plain}{\fancyhead{}}
\title{A review of depth-averaged models for dry granular flows: Savage-Hutter and $\mu(I)$-rheology\thanks{This work was supported by the University of Seville as part of Doc-Course program.}}
\author{FATHI HASSINE\thanks{\textit{UR Analysis and Control of PDE 13ES64, Department of Mathematics, Faculty of Sciences of Monastir, University of Monastir, 5019 Monastir, Tunisia email:} \texttt{fathi.hassine@fsm.rnu.tn}}}
\date{}
\begin{document}
\maketitle
\begin{center}
\abstract{
In this work we present a review of the most popular depth-averaged models to simulate dry granular flows such as aerial avalanches. The classical Savage-Hutter model and recent ones using a $\mu(I)$-rheology law are studied. The objective is firstly to point out the advantages of each such models and secondly to understand how the hypothesis considered in the derivation process influence on the final system.
}
\end{center}
\section{Introduction and physical model}



The increased air pollution of recent years has been responsible for the ecological phenomenon of the dying forests. As a result, the Alpine regions have experienced a greater potential for the occurrence of slow moving landslides, as well as for more rapidly flowing rock, ice and snow avalanches. This has led to an increased emphasis on the study of these and related phenomena.

Direct field observations of the dynamics of rockfalls or avalanches are extremely difficult to make and as a result, all present theoretical models contain certain hypothetical facets. In addition, the field events can be extremely complex in terms of the kinds and sizes of materials that are present as well as the bed and avalanche geometries that might be involved. In view of this complexity and the lack of detailed field data against which theoretical models can be tested, it was decided to approach the problem of avalanche motion in the following way. We have attempted to examine the simplest two-dimensional problem which embodies the main features of avalanche motion from the initiation of motion to the final run out when the material finally comes to rest and to reduce the analysis of this problem to as simple a treatment as possible. In doing so we shall exploit simplicities inherent in the constitutive behavior of flowing granular materials. The plan is to describe the mass dispersion, the velocity field, the run-out distances and deposition areas, and to compare the theoretical predictions with well defined laboratory experiments. A satisfactory fit of such a theory with laboratory data still does not imply that the theory is adequate to describe large scale processes in nature. The laboratory materials used for simulation purposes may be too simple and idealized. There may be large scale effects (high pressure melting of material, etc.) that are present in the field, but which cannot be modeled in the laboratory.

In this notes we are interested to the study of modeling dynamic for areal avalanches in the case where they are considered as dry granular flows. Most of the models devoted to gravitational granular flows describe the behavior of dry granular material  in which a shallow water type model (i.e thin layer approximation for the continuum medium) is derived to describe granular flows over a slopping plane based on Mohr-Coulomb consideration: a Coulomb friction is assumed to reflect the avalanche/bottom interaction and the normal stress tensor is defined by a constitutive law relating the longitudinal and the normal stress through a proportionality factor $K$. One of this kind of models is the Savage-Hutter and the model with $\mu(I)$-rheology that we analyze in this work.

The Savage-Hutter (SH) avalanche model~\cite{SH1,SH2} and its extensions~\cite{GWH,GH,GKH,HG,HKPS,HN,KGH,PH,PEH}, henceforth also called SH-model, is a dynamical fluid-like model which consists of hyperbolic partial differential equations for the distribution of the depth and the (two) topography-parallel, depth averaged velocity components in an avalanching mass of cohesion less granules (e.g., sand, grains, rocks and snow). It is designed to predict the motion and deformation from initiation to run out along a concomitantly determined avalanche track along a prescribed topography. In the past, it has been used to describe flows in straight and curved chutes~\cite{GH, HG, HG,HKPS,HN}, in channels with plane and parabolic cross sections and simply curved thalwegs~\cite{GWH,GKH,H}, but has been extended to flows in corries having arbitrarily curved and twisted thalwegs and arbitrary topographies~\cite{PH,PEH}. The basic simplifying assumptions in the various models are mathematically not exactly the same; however, physically they are identical, namely consisting of the assumption of density preserving (incompressibility), the assumption of shallowness of the avalanche piles, and of small topographic curvatures, the assumptions of Coulomb-type sliding with bed friction angle $\delta$, Mohr-Coulomb frictional behavior in the interior with internal angle of friction $\phi\geq\delta$ and an ad-hoc assumption, reducing the number of Mohr's circles in three dimensional stress states from three to one, and nearly uniform velocity profile through the avalanche depth.

Rheology is the study of the flow of matter, primarily in a liquid state, but also as "soft solids" or solids under conditions in which they respond with plastic flow rather than deforming elastically in response to an applied force. It applies to substances which have a complex micro structure, such as muds, sludges, suspensions, polymers and other glass formers (e.g., silicates), as well as many foods and additives, bodily fluids (e.g. blood) and other biological materials or other materials which belong to the class of soft matter.

Newtonian fluids can be characterized by a single coefficient of viscosity for a specific temperature. Although this viscosity will change with temperature, it does not change with the strain rate. Only a small group of fluids exhibit such constant viscosity. The large class of fluids whose viscosity changes with the strain rate (the relative flow velocity) are called non-Newtonian fluids.

Rheology generally accounts for the behavior of non-Newtonian fluids, by characterizing the minimum number of functions that are needed to relate stresses with rate of change of strain or strain rates. For example, ketchup can have its viscosity reduced by shaking (or other forms of mechanical agitation, where the relative movement of different layers in the material actually causes the reduction in viscosity) but water cannot. Ketchup is a shear thinning material, like yogurt and emulsion paint, exhibiting thixotropy, where an increase in relative flow velocity will cause a reduction in viscosity, for example, by stirring. Some other non-Newtonian materials show the opposite behavior: viscosity going up with relative deformation, which are called shear thickening or dilatant materials.

An other way to introduce the granular behavior in avalanches is to consider the dry granular flow as non-Newtonian fluid involving a definition of the rheology in the system. This is usually made through a nonlinear viscous law in the stress tensor. One of this definition is the $\mu(I)$-rheology. In this work we analyze the difference of these two kind of models: Savage-Hutter and model using the $\mu(I)$-rheology. The objective is firstly is to point out the advantages of each such models and secondly to understand how the hypotheses considered in the derivation of each model influence into the final system.
\section{Derivation process}
\subsection{Starting system}
\begin{figure}
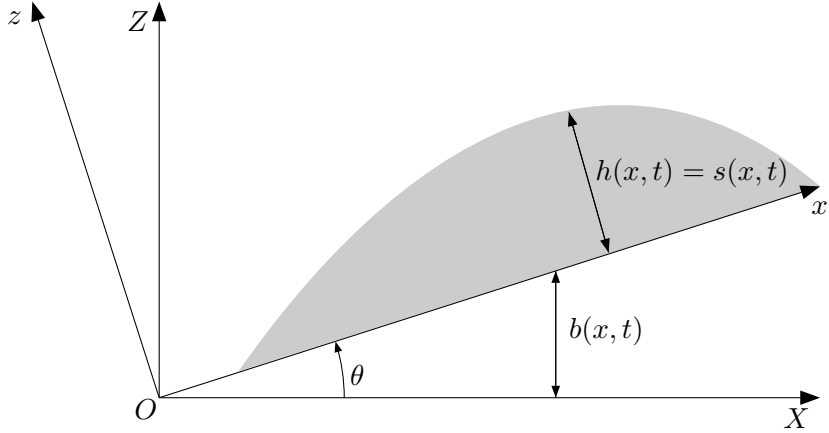

\figinit{pt}
\figpt 0:(0,0)
\figpt 1:(250,0)
\figpt 2:(0,150)
\figpt 3:(250,80)
\figpt 4:(75,0)
\figpt 5:(180,140)
\figpt 6:(100,110)
\figpt 7:(30,10)
\figpt 8:(170,55)
\figpt 9:(155,108)
\figpt 10:(165,85)
\figpt 11:(-48,150)
\figpt 12:(150,0)
\figpt 13:(150,48)
\figpt 14:(155,20)
\psbeginfig{}
\psset arrowhead(fillmode=yes)
\psaxes 0(0,250,0,150)
\psset arrowhead(length=0,ratio=0.2)
\psarrowcircP 0 ; 70 [1,3]
\psset(fillmode=yes,color=0.8)
\psBezier1[3,5,6,7]
\psset(fillmode=no)
\psset(fillmode=no,color=0)
\psset arrowhead(width=0.1,length=0.1,ratio=0.1)
\psarrow [9,8]
\psarrow [8,9]
\psarrow [12,13]
\psarrow [13,12]
\psset arrowhead(width=0.1,length=0.1,ratio=0.05)
\psarrow [0,11]
\psset arrowhead(width=0.1,length=0.1,ratio=0.03)
\psarrow [0,3]
\psendfig
\figvisu{\figBoxA}{}{
\figwritesw 0:$O$(1)
\figwritesw 1:$X$(5.5)
\figwritesw 2:$Z$(5.5)
\figwriten 4:$\theta$(5)
\figwritee 10:$h\textrm{(}x,t\textrm{)}=s\textrm{(}x,t\textrm{)}$(0)
\figwrites 3:$x$(5.5)
\figwritesw 11:$z$(5.5)
\figwritene 14:$b\textrm{(}x,t\textrm{)}$(0)
}
\centerline{\box\figBoxA}
\caption{Local coordinates}
\label{f2}
\end{figure}

Let $Oxz$ be a local rectangular Cartesian coordinate system with the $x$-axis orientated down a slope at an angle $\theta$ to the horizontal and the $z$-axis being the upward pointing normal (see figure~\ref{f2}). The velocity $\textbf{U}$ has components $(u,w)$ in the $\textbf{X}=(x,z)$ direction respectively, and the grains have constant intrinsic density $\rho^{*}$. The solids volume fraction $\Phi$ is assumed to be constant and uniform throughout the material (see~\cite{M}), so the partial density $\rho=\Phi\rho^{*}$ is constant and uniform. The equations are derived from the principles of conservation of mass which implies that the granular material is incompressible and conservation of momentum in two dimension namely,
\begin{equation*}
\begin{split}
\nabla.\textbf{U}=0
\\
\p_{t}(\rho\textbf{U})+\rho\textbf{U}.\nabla_{\textbf{X}}\textbf{U}=\nabla.\sigma+\rho\textbf{g}
\end{split}
\end{equation*}
where $\textbf{g}=(0,-g)$ , being is the gravitational acceleration vector, $\textbf{U}$ is the velocity field and $\rho$ is the density of the granular mass. Moreover, we denote by $\sigma$ the stress tensor
$$
\sigma=\left(\begin{array}{ll}
\sigma_{xx}&\sigma_{xz}
\\
\sigma_{zx}&\sigma_{zz}
\end{array}\right)
$$
with $\sigma_{xz}=\sigma_{zx}$. For Savage-Hutter and $\mu(I)$-rheology models the dependent variables are the fluid height or depth, $h$, and the two-dimensional fluid velocity field, $u$ and $w$. With the proper choice of units, the conserved quantities are mass, which is proportional to $h$, and momentum, which is proportional to $uh$ and $wh$.

The Rheology model is assumed to be with a stress tensor defined with the sum of pressure component and viscous one,
$$
\sigma=-p\mathrm{\textbf{I}}+\tau
$$
where $\mathrm{\textbf{I}}$ is the $2\times2$ identity matrix and the deviatoric part
$$
\tau=\left(\begin{array}{ll}
\tau_{xx}&\tau_{xz}
\\
\tau_{zx}&\tau_{zz}
\end{array}\right).
$$
A fluid is said to be Newtonian if $\tau$ is proportional to the rate deformation tensor $D(\textbf{U})$, where
$$
D(\textbf{U})=\frac{1}{2}\left(\nabla\textbf{U}+{}^{t}\nabla\textbf{U}\right).
$$
Then for a Newtonian fluid, such as water, we have
$$
\tau=\eta D(\textbf{U}),
$$
where $\eta$ is the constant viscosity coefficient, depending on the material.

The behavior of the flows of material like honey, corn flour or paint cannot be modeled with such linear relation. Moreover the concept of Newtonian fluid is an idealization: there are always nonlinear relation between the shear stress and the shear rate. The study of the deviation (from the linear law) of $\tau$ as a function of $D(\textbf{U})$ belongs to the field of Rheology which is defined by the study of deformation and flow of complex fluids.

Fluid in which $\eta$ is a function of $D(\textbf{U})$, $\eta=\eta(D(\textbf{U}))$, i.e.
$$
\tau=\eta(D(\textbf{U})) D(\textbf{U}),
$$
are called Generalization Newtonian fluids. For example corn flour is a material whose behavior is fluid when it is gently mixed but it becomes very viscous if it is strongly mixed. That's viscosity increases with shear, suchmaterails are shear-thinckening. There are many other materials, like paint, whose viscosity decrease with shear. These materials are shear-thinning. Shear-thinckening and shear-thinning can be called using power-law fluids. The viscosity of power-law fluid is defined, for some positive constant $\overline{\eta}>0$ by
$$
\eta(D(\textbf{U}))=\overline{\eta}|D(\textbf{U})|.
$$

The flow of some materials cannot be modeled by a power-law model. This is the case of clay, snow or lava that only flow when the shear stress is bigger than a critical value. These materials are example of what we call threshold fluid. Belong the stress $\tau_{c}$ the material present a rigid behavior but above $\tau_{c}$ the material being to flow. They are visco-plastic materials.

The $\mu(I)$-rheology is a nonlinear viscous law, with a strain-rate invariant and pressure-dependent viscosity, that has proved to be effective at modeling dry granular flows in the intermediate range of the inertial number $I$. The constitutive model for the granular material is provided by the $\mu(I)$-rheology (Jop et al.~\cite{JFP2}and~\cite{JFP1}), which is a nonlinear viscous law with a pressure and strain-rate-dependent viscosity of the form
\begin{equation}\label{1}
\tau=\mu(I)p\frac{D}{\|D\|}
\end{equation}
where $\mu$ is the friction law and $I$ is the inertial number. We remind that the strain-rate $D$ is defined in the term of the velocity gradient as
\begin{equation}\label{2}
D=D(\textbf{U})=\frac{1}{2}\left(\nabla\textbf{U}+{}^{t}\nabla\textbf{U}\right),
\end{equation}
while $\|D\|$ is a second invariant of the strain-rate tensor
\begin{equation}\label{4}
\|D\|=\sqrt{\frac{1}{2}\mathrm{tr}(D^{2})}
\end{equation}
with $\mathrm{tr}$ is the trace. Note that~\eqref{2} is the standard definition of the strain-rate tensor, which differs from the definition used by Jop et al.~\cite{JFP1}, so there is an extra factor of two in the definition of the non-dimensional inertial number
$$
I=\frac{2\|D\|d}{\sqrt{p/\rho^{*}}}
$$
where $d$ is the diameter of the grains. The inertial number is the ratio of the time scale for microscopic rearrangements of the particles at a given confining pressure, $\displaystyle\frac{d}{\sqrt{p/\rho^{*}}}$, to the time scale given by the bulk shear rate $\displaystyle\frac{1}{\|D\|}$ and is equal to the square root of the Savage or Coulomb number (Savage~\cite{S} and Ancey, Coussot and Evesque~\cite{ACE}).

The rate dependence in the rheology~\eqref{1} arises from the increase of the friction coefficient $\mu$ with increasing inertial number $I$. The dependence was determined from basal friction measurements that were made on an inclined plane (Pouliquen~\cite{P1}; Pouliquen and Forterre~\cite{PF}). They observed that steady-uniform depth flows developed between two critical inclination angles $\theta_{1}$ and $\theta_{2}$. For slope angles below $\theta_{1}$ there was no flow and for angles above $\theta_{2}$ the flows accelerated. In the steady-uniform regime they determined the empirical basal friction law
\begin{equation}\label{3}
\mu(Fr,h)=\mu_{1}+\frac{\mu_{2}-\mu_{1}}{\displaystyle\frac{\beta h}{\mathscr{L}Fr}+1},
\end{equation}
where the friction coefficients $\mu_{1}=\tan(\theta_{1})$ and $\mu_{2}=\tan(\theta_{2})$. The parameter $\beta$ is a dimensionless empirical constant (Pouliquen~\cite{P1}), whilst $\mathscr{L}$ has the dimensions of a length and is characterized by the depth of flow over which a transition between the angles $\theta_{1}$ and $\theta_{2}$ occurs and, as such, is dependent on the material properties of the flowing particles and on the bed roughness conditions. On an inclined plane the Froude number
$$
Fr=\frac{\overline{u}}{\sqrt{gh\cos(\theta)}},
$$
is defined as the ratio of the depth-averaged flow velocity $\overline{u}$ to the gravity wave speed $\sqrt{gh\cos(\theta)}$ where $g$ is the constant of gravitational acceleration and $h$ is the flow thickness (see e.g. Gray et al.~\cite{GTN}). In these steady-uniform flows the inertial number, $I$, is constant and there is a Bagnold velocity profile through their depth (see~\cite{M}). Using the fact that the depth-averaged Bagnold velocity is equal to
$$
\overline{u}=\frac{2I}{5d}\sqrt{gh\cos(\theta)}h^{\frac{3}{2}},
$$
Jop et al~\cite{JFP2} substituted for the Froude number and the depth-averaged velocity in~\eqref{3} to obtain a general expression for the friction as a function of the inertial number
$$
\mu(I)=\mu_{1}+\frac{\mu_{2}-\mu_{1}}{I_{0}/I+1},
$$
where the constant
$$
I_{0}=\frac{5\beta d}{2\sqrt{h}\mathscr{L}}.
$$
The basal and internal friction laws are therefore intimately linked. The friction law~\eqref{3} is only strictly valid for Froude numbers above $\beta$. For Froude numbers below this value, Pouliquen and Forterre~\cite{PF} determined a transition law, which plays an important role in the development of static regions (see e.g. Mangeney et al.~\cite{MBTVB} and Johnson and Gray~\cite{JG}).
\subsection{Boundary conditions}
Let $h(x,t)$ be the height of the granular layer along the normal direction to the bed. We denote by $\textbf{n}^{s}$ the unit normal vector to the free granular surface $s(x,t)$ with positive vertical component and $\textbf{n}^{b}$ the unit normal vector to the bottom $b(x,t)$, i.e. by denoting $F^{s}=s(x,t)-z$ and $F^{b}=b(x,t)-z$ then $\textbf{n}^{s}$ and $\textbf{n}^{b}$ are given by
$$
\textbf{n}^{s}=\frac{\nabla F^{s}}{|\nabla F^{s}|}\qquad\text{ and }\qquad \textbf{n}^{b}=\frac{\nabla F^{b}}{|\nabla F^{b}|}.
$$
The granular material is subject to kinematic conditions at its free surface and its base:
\begin{enumerate}
	\item[1/] \textbf{At the free surface:}
	
	For the Savage-Hutter model and the model with $\mu(I)$-rheology we have a kinematic condition at the free surface,
	\begin{equation}\label{SHbc1}
	\p_{t}h^{s}+u^{s}\p_{x}h^{s}-w^{s}=0,
	\end{equation}
	which means that the particles at the free surface are transported with velocity $(u^{s},w^{s})$. And two boundary conditions are imposed,
	\begin{equation}\label{SHbc2}
	\textbf{n}^{s}.\sigma\textbf{n}^{s}=0
	\end{equation}
	\begin{equation}\label{SHbc3}
	\sigma\textbf{n}^{s}-\textbf{n}^{s}(\textbf{n}^{s}.\sigma\textbf{n}^{s})=\left(\begin{array}{c}
	fric_{h}(u)
	\\
	0
	\end{array}\right),
	\end{equation}
	where $fric_{h}(u)$ is the friction term between the granular layer and the air. For the sake of simplicity we will suppose that $fric_{h}{u}=0$.
	
	\item[2/] \textbf{At the base:} Two bondary conditions are imposed for \underline{the Savage-Hutter model}
	\begin{equation}\label{SHbc4}
	\textbf{U}^{b}.n^{b}=0
	\end{equation}
	the condition of non penetration and
	\begin{equation}\label{SHbc5}
	\sigma\textbf{n}^{b}-\textbf{n}^{b}(\textbf{n}^{b}.\sigma\textbf{n}^{b})=\left(\begin{array}{c}
	\displaystyle-\textbf{n}^{b}.\sigma\textbf{n}^{b}\frac{u^{b}}{|u^{b}|}\tan(\delta_{0})
	\\
	0
	\end{array}\right),
	\end{equation}
	the Coulomb friction law, defined in term of the angle of repose $\delta_{0}$ (see~\cite{SH1}).
	
	A no slip condition at the bas is imposed for \underline{the model with $\mu(I)$-rheology}
	\begin{equation}\label{Rbc4}
	\textbf{U}^{b}=0.
	\end{equation}
	The no slip condition is consistent with observations of flows on rough beds made of particles of the same size and shape that are glued to the base (Pouliquen\cite{P1,P2}, Pouliquen and Forterre~\cite{PF} and GDR-MiDi~\cite{M}).
\end{enumerate}
Where here the superscripts `$s$' and `$b$' on the velocity indicate evaluation at the surface and base, respectively.

The main difference is the condition of friction at the bottom, while Savage-Hutter model considers a Coulomb friction law, no slip condition is assumed for the $\mu(I)$-rheology model.
\subsection{Dimensional analysis}
The shallowness of the flow is now exploited in order to obtain simplified depth averaged equations. The avalanche is assumed to be of a typical thickness $H$ which is much smaller than the downslope length scale $L$. This suggests introducing, for both models, non-dimensional variables indicated, by the tide $\tilde{.}$, of the form
\begin{equation*}
x=L\tilde{x},\quad z=H\tilde{z},\quad s=H\tilde{s},\quad h=H\tilde{h},\quad b=H\tilde{b}.
\end{equation*}

Different order has been assumed for downstream flow speeds which are assumed to be of the order $\sqrt{gL}$, time parameter and the Cauchy stress tensor parameters which are chosen for the Savage-Hutter model as follow
\begin{equation*}
\begin{array}{c}
\displaystyle t=\sqrt{\frac{L}{g}}\tilde{t},
\\
u=\sqrt{gL}\tilde{u},\quad w=\epsilon\sqrt{gL}\tilde{w},
\\
\sigma_{xx}=gH\tilde{\sigma}_{xx},\quad\mathscr{\sigma}_{zz}=gH\tilde{\sigma}_{zz},\quad \sigma_{xz}=\sigma_{zx}=gH\mu\tilde{\sigma}_{xz}=gH\mu\tilde{\sigma}_{zx},
\end{array}
\end{equation*}
where $\mu=\tan(\delta_{0})$, $\delta_{0}$ being the angle of repose in the Coulomb term and the aspect ratio $\epsilon= H/L$, which is supposed to be small.

For the model with the $\mu(I)$-rheology typical downstream flow speeds are assumed to be of the order of the gravity wave speed $\sqrt{gH}$, and mass balance implies that typical normal velocities in the $z$ direction are of magnitude $\epsilon\sqrt{gH}$. The pressure scaling $\rho gh$ is based on a lithostatic balance in the normal momentum equation
\begin{equation*}
u=\sqrt{gH}\tilde{u},\quad w=\epsilon\sqrt{gH}\tilde{w},
\end{equation*}
Typical magnitudes for the strain-rate and, hence, the deviatoric stresses can then be determined from the constitutive relations~\eqref{1}-\eqref{4}. This suggests introducing non-dimensional variables, indicated by the hat, of the form
\begin{equation*}
\begin{array}{c}
\displaystyle t=\frac{L}{\sqrt{gH}}\tilde{t},\quad p=\rho g H\tilde{p},
\\
\displaystyle D_{xx}=\epsilon\sqrt{\frac{g}{H}}\tilde{D}_{xx},\; D_{zz}=\epsilon\sqrt{\frac{g}{H}}\tilde{D}_{zz},\; D_{xz}=D_{zx}=\sqrt{\frac{g}{H}}\tilde{D}_{xz}=\sqrt{\frac{g}{H}}\tilde{D}_{zx},
\\
\displaystyle\tau_{xx}=\epsilon\rho gH\tilde{\tau}_{xx},\; \tau_{zz}=\epsilon\rho gH\tilde{\tau}_{zz},\; \tau_{xz}=\tau_{zx}=\rho gH\tilde{\tau}_{xz}=\rho gH\tilde{\tau}_{zx},\;\sigma=\rho gH\tilde{\sigma}.
\end{array}
\end{equation*}
\subsection{Assumptions}
The assumptions made to the Savage-Hutter model are
\begin{itemize}
	\item The following constitutive law is considered(see~\cite{SH1})
	\begin{equation}\label{SHa1}
	\sigma_{xx}=K\sigma_{zz},
	\end{equation}
	where $K$ measures the entropy or normal stress effect: where $K=1$ corresponds to isotropic conditions, $K\neq 1$ makes "overburden pressures" different from the normal stresses parallel to the basical surface. In this case the Shallow Water equation $K=1$ is assumed. The coefficient $K$ is defined according to the motion of the granular layer:
	$$
	K=\left\{\begin{array}{ll}
	K_{act}&\text{if }\displaystyle\frac{\p u}{\p x}>0,
	\\
	K_{pas}&\text{if }\displaystyle\frac{\p u}{\p x}<0,
	\end{array}\right.
	$$
	where
	$$
	K_{act/pas}=\frac{2}{\cos^{2}(\phi)}\left(1\mp\left(1-\frac{\cos^{2}(\phi)}{\cos^{2}(\delta_{0})}\right)\right)-1
	$$
	being $\phi$ the internal friction angle, defined in terms of the type of grains and size.
	\item In~\cite{G} Gray introduced the assumption that the Coulomb term is of order $\gamma$ for some $\gamma\in(0,1)$. That is,
	\begin{equation}\label{SHa2}
	\mu=\tan(\delta_{0})=\mathscr{O}(\epsilon^{\gamma}).
	\end{equation}
	\item Supposing finally a constant profile of the velocities.
\end{itemize}
The assumptions made to the model with $\mu(I)$-rheology are:
\begin{itemize}
	\item The sign of $\displaystyle\frac{\p\overline{u}}{\p\tilde{z}}$ is constant, where we denoted by $\overline{u}$ the average of the velocity along the normal direction:
	$$
	\overline{u}=\frac{1}{h}\int_{b}^{s}u(x,z)\,\ud z.
	$$
	\item The gravitational force $\tilde{h}\sin(\theta)$ and the basal friction $\mu\tilde{h}\cos(\theta)$ are both order-unity quantities, their difference is typically much smaller. To formalize this, it is assumed that
	\begin{equation}\label{Ra1}
	\tilde{h}\sin(\theta)-\mu\tilde{h}\cos(\theta)\mathrm{sgn}(\tilde{\overline{u}})=\epsilon\tilde{h}\cos(\theta)\left(\tan(\theta)-\mu(Fr,\tilde{h})\mathrm{sgn}(\tilde{\overline{u}})\right)+\mathscr{O}(\epsilon^{2}),
	\end{equation}
	i.e. gravity balances friction to leading order and their difference is small.
\end{itemize}
\section{Final models}
The final system of equations of the Savage-Hutter model is done as follows
\begin{equation}\label{SHsys}
\left\{\begin{array}{c}
\displaystyle\frac{\p h}{\p t}+\frac{\p}{\p x}(h\overline{u})=0,
\\
\displaystyle\frac{\p}{\p t}(h\overline{u})+\frac{\p}{\p x}\left(h\overline{u}^{2}+\frac{1}{2}g\cos(\theta)h^{2}K\right)=\mathscr{S}_{1},
\end{array}\right.
\end{equation}
where 
$$
\mathscr{S}_{1}=-gh\left(\cos(\theta)\frac{\p b}{\p x}+\sin(\theta)\right)+\mathscr{T}_{1}
$$
and $\mathscr{T}_{1}$ represents the Coulomb friction term. This term must be understood as follows
\begin{equation*}
\begin{array}{lll}
\text{If }|\mathscr{T}_{1}|\geq\sigma_{c}&\Longrightarrow& \displaystyle\mathscr{T}_{1}=-gh\cos(\theta)\frac{\overline{u}}{|\overline{u}|}\tan(\delta_{0}),
\\
\text{If }|\mathscr{T}|<\sigma_{c}&\Longrightarrow& U=0,
\end{array}
\end{equation*}
where $\sigma_{c}=gh\cos(\theta)\tan(\delta_{0})$, and which originates from the Cauchy stress tensor $\sigma$ combined with the boundary conditions~\eqref{SHbc2}-\eqref{SHbc3}, ~\eqref{SHbc4}-\eqref{SHbc5}, the constitutive law~\eqref{SHa1} and the Gray assumption~\eqref{SHa2}. 

The final system of equations of the model of $\mu(I)$-rheology is done as follows
\begin{equation}\label{Rsys}
\left\{\begin{array}{c}
\displaystyle\frac{\p h}{\p t}+\frac{\p}{\p x}(h\overline{u})=0,
\\
\displaystyle\frac{\p}{\p t}(h\overline{u})+\frac{\p}{\p x}\left(\chi h\overline{u}^{2}+\frac{1}{2}g\cos(\theta)h^{2}\right)=\mathscr{S}_{2}+\frac{\p}{\p x}\left(\nu h^{\frac{3}{2}}\frac{\p\overline{u}}{\p x}\right),
\end{array}\right.
\end{equation}
where the shape factor $\chi$ is the ratio of the depth-averaged square of the velocity which is formally equal to $\displaystyle\frac{5}{4}$ for the Bagnold velocity profile, but in virtually all granular flow models $\chi$ is assumed to be unity for simplicity. The dimensional source term is
\begin{equation}\label{5}
\mathscr{S}_{2}=-gh\left(\cos(\theta)\frac{\p b}{\p x}-\sin(\theta)\right)+\mathscr{T}_{2},
\end{equation}
which is the combination of gravity acceleration, effective basal friction and topography gradients, where
$$
\mathscr{T}_{2}=-gh\mu(Fr,h)\mathrm{sgn}(\overline{u})\cos(\theta).
$$
The dimensional coefficient in the viscous law
$$
\nu=\frac{2}{9}\frac{\mathscr{L}\sqrt{g}}{\beta}\frac{\sin(\theta)}{\sqrt{\cos(\theta)}}\left(\frac{\tan(\theta_{2})-\tan(\theta)}{\tan(\theta)-\tan(\theta_{1})}\right).
$$

Both equations~\eqref{SHsys} and~\eqref{Rsys} are the familiar shallow-water-type avalanche equations (resulting from the surface kinematic conditions~\eqref{SHbc1}, downslope surface traction condition~\eqref{SHbc2}-\eqref{SHbc3} or the no penetration condition~\eqref{SHbc4} and the Coulomb friction law), which are commonly used in the literature and have proved their effectiveness over many years (e.g. Grigorian et al. \cite{GEI} Savage and Hutter \cite{SH2}, Gray et al. \cite{GA}, Pouliquen \cite{P2}, Pouliquen and Forterre~\cite{PF}, Gray et al.~\cite{GTN}). The main difference  between the classical Shallow-Water equations and the Savage-Hutter equation model is the presence of the Coulomb friction term: if a closed domain is considered and the Coulomb friction term is neglected, the stationary solution is horizontal free surface, corresponding to water at rest. 
It is interesting how, under the assumption~\eqref{Ra1}, the combination of the no slip condition~\eqref{Rbc4} and the internal rheology~\eqref{1} naturally give rise to an effective basal friction in the source terms~\eqref{5}. This is the only effect of the rheology on the flow, as the depth-averaged in-plane deviatoric stress gradient does not contribute to the leading-order momentum balance.

For the large majority of situations, the new depth-averaged $\mu(I)$-rheology can be neglected, but when sharp gradients in $\overline{u}$ develop, a boundary-layer forms in which the viscous terms play a significant role. This system therefore has all the advantages of the classic shallow-water-type avalanche models (e.g. Grigorian et al. \cite{GEI}, Savage and Hutter \cite{SH2}, Gray et al. \cite{GWH} Pouliquen \cite{P2}; Pouliquen and Forterre~\cite{PF}; Gray et al.~\cite{GTN}), but has the extra physics necessary to obtain highwavenumber cutoff (Forterre~\cite{F}) as well as the ability to regularize ill-posed models (Woodhouse et al.~\cite{WTJKG}). The term $\displaystyle\frac{\nu}{2} h^{3/2}$ is the coefficient of depth-averaged viscosity. The dependence on the thickness to the three halves power is a direct result of the $\mu(I)$-rheology~\eqref{1}.

The main difference of the two models is 
\begin{itemize}
	\item firstly the presence of the viscous term  in the $\mu(I)$-rheology equation due to the definition of the Cauchy stess tensor and which is a direct result of the $\mu(I)$-rheology,
	\item and secondly the source terms $\mathscr{S}_{1}$ and $\mathscr{S}_{2}$ whose terms are the same except that instead of having $\tan(\delta_{0})$ in $\mathscr{S}_{2}$ as in $\mathscr{S}_{1}$ we have $\mu(I)$.
\end{itemize} 

The advantages of the $\mu(I)$-rheology model are
\begin{itemize}
	\item First, we have more properties of the dry granular material,
	\item Second, the definition of $\mathscr{T}_{2}$ is given, contrary to that of $\mathscr{T}_{1}$ which is non closed.
\end{itemize}
The main disadvantage of the $\mu(I)$-rheology model is the assumption that $\|D(u)\|\neq 0$, which is not necessarily true.

\subsubsection*{Acknowledgments}
Sincere thanks to professor Enrique D. Fernandez Nieto and professor Gladys Narbona Reina for inspiring question, their greatly contribution to this work and for their careful reading of the manuscript. I want to thank also all the organizing committee of the Doc-course for their hospitality and for giving me this great opportunity to spend an amazing stay in Spain and expressing myself in differently.
\nocite{*}
\bibliographystyle{alpha}
\bibliography{BibCompMod}
\addcontentsline{toc}{section}{References}
\end{document}